\documentclass[VANCOUVER,STIX1COL]{WileyNJD-v2}

\usepackage{amsfonts}
\usepackage{amsmath}
\usepackage{graphicx}
\usepackage{subfigure}
\usepackage{bm}         \usepackage{multirow}   \usepackage{siunitx}    \usepackage{soul}       
\usepackage[capitalise,noabbrev]{cleveref}
\crefformat{equation}{(#2#1#3)}

\usepackage{tikz}
\usepackage[most]{tcolorbox}
\usetikzlibrary{plotmarks}

\renewcommand{\vec}[1]{\bm{#1}}

\newtcbox{\highlight}[2]{enhanced, box align=base, nobeforeafter, colback=#1, colframe=#2, size=small, left=0pt, right=0pt, boxsep=2pt}

\newcommand{\qtwo}{$\mathbb{Q}_2$}
\newcommand{\isoqtwo}{$\mathbb{Q}_1 \text{iso}\kern1pt\mathbb{Q}_2$}
\newcommand{\qtwoqone}{$\boldsymbol{ \mathbb{Q}}_2/\mathbb{Q}_1$}

\newcommand{\isoqtwoqone}{$\boldsymbol{ \mathbb{Q}}_1 \text{iso}\kern1pt\boldsymbol{ \mathbb{Q}}_2/ \mathbb{Q}_1$}
\newcommand{\isoptwopone}{$\boldsymbol{ \mathbb{P}}_1 \text{iso}\kern1pt\boldsymbol{ \mathbb{P}}_2/ \mathbb{P}_1$}

\articletype{}

\begin{document}

\newcommand{\thistitle}{~Low-order preconditioning of the Stokes equations}
\title{\thistitle}

\author[1]{Alexey Voronin}
\author[2]{Yunhui He}
\author[3]{Scott MacLachlan}
\author[1]{Luke N. Olson}
\author[4]{Raymond Tuminaro}

\authormark{Voronin, He, MacLachlan, Olson, Tuminaro}

\address[1]{\orgdiv{Department of Computer Science}, \orgname{University of Illinois at Urbana-Champaign}, \orgaddress{\city{Urbana} \state{IL}, \country{USA}}}
\address[2]{\orgdiv{Department of Applied Mathematics}, \orgname{University of Waterloo}, \orgaddress{\city{Waterloo} \state{ON}, \country{Canada}}}
\address[3]{\orgdiv{Department of Mathematics and Statistics}, \orgname{Memorial University of Newfoundland}, \orgaddress{\city{St.~John's} \state{NL}, \country{Canada}}}
\address[4]{\orgdiv{Computational Mathematics}, \orgname{Sandia National Laboratories}, \orgaddress{\city{Livermore} \state{CA}, \country{USA}}}
\corres{*Alexey Voronin. Department of Computer Science.  University of Illinois at Urbana-Champaign. \email{voronin2@illinois.edu}}

\abstract[Summary]{A well-known strategy for building effective preconditioners for higher-order discretizations of some PDEs, such as Poisson's equation, is to leverage effective preconditioners for their low-order analogs. In this work, we show that high-quality preconditioners can also be derived for the Taylor-Hood discretization of the Stokes equations in much the same manner. In particular, we investigate the use of geometric multigrid based on the \isoqtwoqone{} discretization of the Stokes operator as a preconditioner for the \qtwoqone{} discretization of the Stokes system. We utilize local Fourier analysis to optimize the damping parameters for Vanka and Braess-Sarazin relaxation schemes and to achieve robust convergence. These results are then verified and compared against the measured multigrid performance. While geometric multigrid can be applied directly to the \qtwoqone{} system, our ultimate motivation is to apply algebraic multigrid within solvers for \qtwoqone{} systems via the \isoqtwoqone{} discretization, which will be considered in a companion paper.}

\keywords{Monolithic Multigrid, Stokes Equations, Braess-Sarazin, Additive Vanka, Local Fourier Analysis}

\maketitle

\section{Introduction}\label{sec:intro}

This paper focuses on developing efficient algorithms for the numerical approximation of solutions to the Stokes equations, which are used to simulate incompressible viscous flow and whose discretization results in saddle-point linear systems~\cite{elman2014finite}. Linear systems of saddle-point type appear in a variety of scientific and engineering applications~\cite{MBenzi_GHGolub_JLiesen_2005a}.  The block structure and indefiniteness of these systems often make it particularly challenging to construct efficient numerical schemes. The coupled physical fields, pressure and velocity, are often discretized using staggered grids or methods where degrees of freedom are not located at the same spatial point, further complicating the development of solvers such as monolithic multigrid for these systems.

Monolithic multigrid methods, that apply multigrid to the entire system in a coupled (or all-at-once) manner, have long demonstrated robust convergence for Stokes problems~\cite{ABrandt_NDinar_1979a, ABrandt_1984a, SPVanka_1986a, JLinden_etal_1989a, ANiestegge_KWitsch_1990a, Braess, VJohn_LTobiska_2000a, MWabro_2004a, MWabro_2006a, MLarin_AReusken_2008a, AJanka_2008a, BGmeiner_etal_2016a, adler2017preconditioning, prokopenko2017algebraic, PFarrelletal2019a}. Most of these approaches are geometric multigrid (GMG) methods, where the multigrid hierarchy is composed of a sequence of coarser discretizations on nested meshes, connected by canonical interpolation operators for each field. There are many fewer approaches to monolithic algebraic multigrid (AMG) for Stokes and similar systems~\cite{MWabro_2004a, MWabro_2006a, AJanka_2008a, prokopenko2017algebraic}, as algebraic coarsening of these matrices is complicated by the presence of higher-order finite-element bases, which are present in a majority of stable Stokes discretizations. This is in part due to two main obstacles. First, the stiffness matrices constructed with higher-order basis functions are no longer M-matrices, for which AMG methods were originally intended~\cite{osti_1660522, heys2005algebraic}. Secondly, the coarsening of one field cannot be done independently of the other, due to \textit{inf-sup} stability concerns~\cite{MWabro_2004a, MWabro_2006a, AJanka_2008a, prokopenko2017algebraic}. In order to address the first obstacle, we propose to coarsen with respect to the element order ($p$) and, then, to employ geometric multigrid on the low-order system by coarsening spatially in $h$.  This does require that a low-order system be constructed, however, it is facilitated by the standard building blocks needed to build the high-order discretization.  While the approach presented here is one of geometric multigrid, we propose this primarily as a building block that can be leveraged in future AMG work that will aim to address the second obstacle.

Low-order preconditioning is a common approach when developing preconditioners for higher-order discretizations. These methods have been successfully used for preconditioning discretized Laplace operators in both geometric and algebraic multigrid contexts~\cite{deville1990finite, heys2005algebraic, napov2014algebraic, olson2007algebraic, orszag1979spectral, xu1996auxiliary}. In the case of the Stokes problem, non-nested low-order geometric multigrid methods have been shown to be effective preconditioners for higher-order discretizations~\cite{john2002non}. In this paper, we construct and use a stable lower-order \isoqtwoqone{} discretization within a monolithic multigrid preconditioner for the \qtwoqone{} discretized Stokes system. While ``double discretization'' or ``defect correction'' schemes~\cite{WHackbusch_1982a, BKoren_1990a} and multigrid methods based on low-order discretizations~\cite{MWabro_2004a, benzi2006augmented, emamilehrstuhl} have been considered before independently of each other, to our knowledge, these schemes have not yet been considered for the stable Taylor-Hood finite-element discretization. The goal of our work is to gauge the effectiveness of this type of preconditioning for this discretization.

A major complication in defect-correction algorithms is the plethora of algorithmic parameters that can arise, including relaxation parameters for both higher-order and lower-order discretizations, and additional parameters related to cycling strategies between the discretizations.  Here, we consider both the choice of and sensitivity to these parameters in the context of geometric multigrid, which we use to approximate the inverse of the \isoqtwoqone{} system within a standard defect-correction scheme.  Using geometric multigrid offers an advantage in that parameter choices can be determined using local Fourier analysis~\cite{Trottenberg2001, wienands2004practical}.  Local Fourier analysis has been studied and applied to many interesting problems, including monolithic multigrid methods for the Stokes equations.  Indeed, early work in this direction focused on the use of both distributed~\cite{ANiestegge_KWitsch_1990a} and Vanka~\cite{Siv_1991a} relaxation for the staggered marker-and-cell (MAC) finite-difference discretization scheme for the Stokes equations.  More recent Fourier analysis has included that of multiplicative~\cite{SPMacLachlan_CWOosterlee_2011a} and additive~\cite{PFarrelletal2019a} Vanka relaxation for the \qtwoqone{} discretization, as well as for Braess-Sarazin relaxation~\cite{he2019local}.  Here, we make use of these tools to optimize multigrid parameters in stationary iterations for the proposed defect-correction scheme.

The remainder of the paper is structured in the following manner. In~\cref{sec:disc}, we briefly introduce the \qtwoqone{} and \isoqtwoqone{} discretizations for the Stokes equations. \Cref{sec:multigrid} describes the monolithic multigrid framework with Braess-Sarazin and additive Vanka relaxation for the Stokes equations, while \Cref{sec:LFA} introduces local Fourier analysis (LFA) as a tool to identify optimal GMG parameters. \Cref{sec:num} exhibits the main contributions of this paper: optimized LFA two-grid convergence factors and the measured multigrid convergence for the Stokes problem.  \Cref{sec:conclusion} presents the conclusions.  The main contribution of this work is that of demonstrating that the principle of low-order preconditioning is also applicable to a mixed finite-element discretization such as the \qtwoqone{} discretization of Stokes.  In our view, this opens a new line of research into the development of truly algebraic multigrid methods for such discretizations, by providing a route to effectively precondition higher-order discretizations without directly applying AMG to them.  A central question in that direction that we aim to address in future work is how to algebraically create a hierarchy of stable coarse-grid operators that allow effective multigrid cycling on problems such as the \isoqtwoqone{} discretization.

\section{Discretization and Solutions of Stokes Equations}\label{sec:disc}

In this paper, we consider the two-dimensional Stokes equations, given by
\begin{subequations}\label{equation:stokes-eq}
\begin{align}
        - \nabla^2 \vec{u} +\nabla p &=\vec{f}  \\
        -\nabla \cdot \vec{u} &= 0
\end{align}
\end{subequations}
where $\vec{u}$ is the velocity of the fluid, $p$ is the pressure, and $\vec{f}$ is a forcing term.  We assume homogenous Dirichlet conditions on $\vec{u}$ over the boundary of a domain $\Omega$  for simplicity.
In addition, we consider a uniform grid of size $h$ over $\Omega$ and finite dimensional spaces of the form
$\vec{\mathcal{X}}^h \subset \vec{H}^1_0(\Omega)$ and $\mathcal{M}^h \subset L_2(\Omega)$, where $\vec{\mathcal{X}}^h$ satisfies the appropriate homogeneous Dirichlet boundary conditions.  The resulting discrete weak formulation of~\cref{equation:stokes-eq} is to
find  $\vec{u} \in \vec{\mathcal{X}}^h$ and $p \in \mathcal{M}^h$ such that
\begin{subequations}\label{equation:stokes-disc}
  \begin{alignat}{2}
    & a(\vec{u}, \vec{v})+b(p, \vec{v}) &&= F(\vec{v})\\
    & b(q, \vec{u})                     &&= 0,
  \end{alignat}
\end{subequations}
for all $q \in \mathcal{M}^h$ and $\vec{v} \in \vec{\mathcal{X}}^h$.  Here, $a(\cdot,\cdot)$ and $b(\cdot,\cdot)$ are bilinear forms and $F(\cdot)$ is a linear form given by
\begin{equation*}
        a(\vec{u}, \vec{v}) = \int_{\Omega} \nabla \vec{u} : \nabla \vec{v},\quad
        b(p, \vec{v}) = - \int_{\Omega} p \nabla \cdot \vec{v},\quad
        F(\vec{v}) = \int_{\Omega} \vec{f} \cdot \vec{v}.
\end{equation*}
An \textit{inf-sup} condition on the finite-dimensional spaces  $\vec{\mathcal{X}}^h$ and $\mathcal{M}^h$ is sufficient to guarantee the uniqueness of the solution up to a constant pressure~\cite{elman2014finite}.

Our focus is on two types of stable mixed finite-element discretizations for $\vec{\mathcal{X}}^h$ and $\mathcal{M}^h$.
The first is the \qtwoqone{} discretization (also known as the \textit{Taylor-Hood} discretization), which uses a bilinear representation of the pressure and a biquadratic representation for the velocity components. The second is the \isoqtwoqone{} discretization, which replaces the $\boldsymbol{\mathbb{Q}}_2$ space for velocities with an $\boldsymbol{\mathbb{Q}}_1$ approximation on a once-refined mesh. This is obtained by overlaying the higher-order nodes with a lower-order mesh~---~see~\cref{fig:meshes}. The \isoqtwoqone{} discretization is not well-known due to its relatively low accuracy and computational efficiency.  Yet, the pairing is known to be \textit{inf-sup} stable~\cite{ern2004book}, and we argue that it can be highly effective when used within a preconditioner for the \qtwoqone{} discretization, by leveraging efficient first-order multigrid solvers.
\begin{figure}[!ht]
   \centering
   \newcommand{\mygrid}[3]{   \def\n{#1}
   \def\m{#2}
   \foreach \y in {0,2,...,\n}{      \draw[black, thick] (0, \y) -- (\n, \y);
   }
   \foreach \x in {0,2,...,\n}{      \draw[black, thick] (\x, 0) -- (\x, \n);
   }

   \ifnum#3=1,{      \foreach \y in {1,2,...,\n}{         \draw[black, densely dotted] (0, \y) -- (\n, \y);
      }
      \foreach \x in {1,2,...,\n}{         \draw[black, densely dotted] (\x, 0) -- (\x, \n);
      }
   }
   \else{}\fi

   \def\h{0.15}
   \foreach \y in {0,2,...,\n}{      \foreach \x in {0,2,...,\n}{         \draw[thick, draw=red!50!red,fill=white] (\x-\h, \y-\h) rectangle +(2*\h,2*\h);
      }
   }
   \foreach \y in {0,1,...,\n}{      \foreach \x in {0,1,...,\n}{
         \pgfmathtruncatemacro\resultX{Mod(\x,2)==0?0:1}
         \pgfmathtruncatemacro\resultY{Mod(\y,2)==0?0:1}
         \pgfmathtruncatemacro\resultXY{\resultX*\resultY}

         \ifnum\resultX=0,{            \ifnum\resultY=0,{               \draw[thick, draw=black!60, fill=black!60] (\x, \y) circle (2pt);
            }
            \else{              \node[mark size=2.pt,color=black!100, thick, fill=white] at (\x,\y) {\pgfuseplotmark{o}};
            }\fi
      }
      \else{         \ifnum\resultXY=1,{            \node[mark size=3pt,color=black!60] at (\x,\y) {\pgfuseplotmark{diamond*}};
         }
         \else{            \ifnum\resultX=1,{               \node[mark size=3pt,color=black!60] at (\x,\y) {\pgfuseplotmark{triangle*}};
            } \else {}\fi
         }\fi
      }\fi
      }
    }
  }

  \subfigure[\qtwoqone]{   \begin{tikzpicture}[x=20pt,y=20pt]
      \begin{scope}[shift={(10pt, 10pt)}]
         \mygrid{6}{2}{0}
      \end{scope}
   \end{tikzpicture}
   }
   \hspace{20mm}
   \subfigure[\isoqtwoqone]{   \begin{tikzpicture}[x=20pt,y=20pt]
      \begin{scope}[shift={(10pt, 10pt)}]
         \mygrid{6}{2}{1}
      \end{scope}
   \end{tikzpicture}
   }
   \newcommand{\mytriangle}{\tikz[baseline=-0.5ex]{\node[mark size=3pt,color=black!60] at (0,0) {\pgfuseplotmark{triangle*}};}}
   \newcommand{\mydiamond} {\tikz[baseline=-0.5ex]{\node[mark size=3pt,color=black!60] at (0,0) {\pgfuseplotmark{diamond*}};}}
   \newcommand{\myasterisk}{\tikz[baseline=-0.5ex]{\node[mark size=2.pt, thick] at (0,0) {\pgfuseplotmark{o}};}}
   \newcommand{\mycircle}  {\tikz[baseline=-0.5ex]{\draw[thick, draw=black!60, fill=black!60] (0, 0) circle (2pt);}}
   \caption{Meshes and degrees of freedom for the \qtwoqone{} and \isoqtwoqone{} discretizations. Dark markers [\protect\mycircle, \protect\myasterisk, \protect\mytriangle, \protect\mydiamond] correspond to the velocity locations (two components of velocity per marker), and red squares \protect\tikz[baseline=-0.5ex]{\protect\draw[thick, draw=red!50!red,fill=white] (-0.10, -0.10) rectangle +(.2,.2);} correspond to the pressure locations.}\label{fig:meshes}
\end{figure}

The challenge presented by the discretization in~\cref{equation:stokes-disc} is the solution of the resulting saddle-point system of the form
\begin{equation}\label{equation:saddle}
    K
    \begin{bmatrix}
        \vec{u}\\
        p
    \end{bmatrix}
    =
    \begin{bmatrix}
        A  & B^T \\
        B  &  0
    \end{bmatrix}
    \begin{bmatrix}
        \vec{u}\\
        p
    \end{bmatrix}
    =
    \begin{bmatrix}
        \vec{f}\\
        0
    \end{bmatrix}
    =b,
\end{equation}
where matrix $A$ corresponds to the discrete vector-Laplacian, and $B$ represents the negative of the discrete divergence operator.  Here, we overload the notation and use $\vec{u}$ and $p$ to denote the discrete velocities and pressure for the remainder of the paper.  We note that, for the boundary conditions described above, the matrix $K$ is singular, with a one-dimensional nullspace associated with constant shifts in the pressure.  While many ``fixes'' to this nullspace are used in practice (e.g., fixing the value of the pressure at a point, or explicitly adding a constraint that the integral of the pressure is zero), we will leave the system in its singular form and implicitly treat the nullspace in our multigrid method.

\section{Monolithic Multigrid}\label{sec:multigrid}

We study a monolithic geometric multigrid approach with standard coarsening~---~i.e., doubling $h$ on each level of the hierarchy.  In the following, we assess the error reduction in multigrid when using a coupled relaxation method followed by a coarse-grid correction.  Interpolation and restriction are constructed using the finite element basis.  In the following, we use $\ell$ to denote the level of the multigrid hierarchy with $\ell=0$ denoting the finest-level problem.  The subscript may be dropped when context is clear.

For a linear system $K_{\ell} u_{\ell}=b_{\ell}$ on level $\ell$ of the multigrid hierarchy, we express the relaxation scheme as a fixed-point iteration
\begin{equation*}\label{eq:fixed_point}
    u_{\ell}^{j+1} = (I-\omega_{\ell} M_{\ell}^{-1} K_{\ell}) u_{\ell}^j + \omega_{\ell}M_{\ell}^{-1}b_{\ell},
\end{equation*}
where $M_{\ell}$ is an inexpensive approximation to $K_{\ell}$ so that the action of its inverse is easy to apply. Here, $(I-\omega_\ell M_{\ell}^{-1} K_{\ell})$ is the associated error-propagation operator for relaxation on $K_{\ell}$, where $\omega_\ell$ is a damping parameter. Using restriction $R_{\ell}$ and interpolation $P_{\ell}$, the coarse-grid correction operator projects the error onto grid $\ell+1$, where a correction is computed. Combining pre- and post-relaxation together with coarse grid correction we arrive at the two-grid error-propagation operator
\begin{equation}\label{eq:TG-operator}
 G_{\ell} = (I - \omega_{\ell} M_{\ell}^{-1} K_{\ell}) (I - P_{\ell} K_{\ell+1}^{-1} R_{\ell} K_{\ell}) ( I -\omega_{\ell} M_{\ell}^{-1} K_{\ell})
\end{equation}
when one sweep of pre- and one sweep of post-relaxation is applied within the multigrid cycle.
The coarse-grid operators for the \isoqtwoqone{} system are computed via the Galerkin
product $K_{{\ell}+1}=R_{\ell} K_{\ell} P_{\ell}$, which is equivalent to rediscretization on the coarse grids.
The correction on the coarsest grid is computed using a pseudoinverse, to account for the pressure nullspace.

\subsection{Low-order Preconditioner}

Lower-order preconditioning has been successfully used for a variety of problems, ranging from elliptic PDEs to saddle-point systems~\cite{deville2002high,heys2005algebraic,orszag1979spectral,olson2007algebraic}. For example, taking $A_{0}$ to be the matrix from a high-order spectral discretization of the Laplace operator and $A_1$ to be the matrix from a second-order finite-difference discretization on an auxiliary mesh constructed       from the nodal points of the high-order mesh, then the condition number of $A_1^{-1} A_0$ is bounded asymptotically by $\pi^2/4$~\cite{deville1990finite, deville1992fourier}, establishing $A_1$ as an effective preconditioner for $A_0$.

Here, we consider a multigrid defect-correction method that incorporates this notion of high- and low-order operators.  This is constructed by taking $K_0$ and $K_1$ to be the higher-order and lower-order Stokes operators assembled using the \qtwoqone{}~and \isoqtwoqone~finite element spaces, respectively.  Note that the \isoqtwoqone{} discretization is not a Galerkin coarsening of the \qtwoqone{} discretization, as the finite-element spaces are not nested.
That is, $K_{1} \neq R_{0} K_{0} P_{0}$ within this multigrid hierarchy.
Since the unknowns for levels 0 and 1 are, however, co-located, we transfer residuals and corrections between the low- and high-order discretizations using identity operators, effectively taking $P_0 = I$ and $R_0 =I$. On level 1, we then use $\gamma$ iterations of monolithic GMG based on the $K_1$ system, which is denoted as $G_1$.  The level-1 multigrid cycle is represented by $\tau(I-G^{\gamma}_{1})K_1^{-1}$, which leads to a generalization of~\cref{eq:TG-operator} for $\ell=0$
\begin{equation}\label{eq:TG-operator-full}
    E=\left(I - \omega_0 M_0^{-1} K_0\right)^{\nu_2}
      \left(I - \tau(I-G_1^\gamma) K_1^{-1} K_0\right)
      \left(I - \omega_0 M_0^{-1} K_0\right)^{\nu_1},
\end{equation}
 where $M_0$ is a relaxation operator based on the $K_0$ system, $\omega_0$ is the relaxation parameter,  $\gamma$ is the number of level-1 $h$-multigrid cycles used to approximately solve the $K_1$ problem,  $\tau$  is a damping parameter for the coarse correction, and  the exponents $\nu_1$ and $\nu_2$ are the number of pre- and post- relaxation sweeps.

The multigrid scheme in~\cref{eq:TG-operator-full} is depicted in~\cref{fig:multilevel}.
If $\gamma=0$, then the method relaxes \textit{only} on the \qtwoqone{} problem and the \isoqtwoqone{} discretization is unused.  At the same time, if $G_1^{\gamma}\equiv 0$ in~\cref{eq:TG-operator-full}, then the \isoqtwoqone{} problem is solved \textit{exactly}.
The left image in \cref{fig:multilevel} illustrates the case where
the overall multigrid cycle employs two-level (blue) and multilevel (green) $h$-multigrid schemes for the \isoqtwoqone{} part of the cycle with $\gamma=1$. The right image shows $\gamma=2$ for a two-level scheme.  In the following tests, we study  multigrid convergence over a range of parameters, possibly omitting relaxation on different levels. For the remainder of the paper, the terms two-level, and multilevel always refer to the number of levels associated with only the GMG (or $h$-multigrid) part of the overall cycle defined by~\cref{eq:TG-operator-full}.
\begin{figure}[!ht]
  \centering
  \definecolor{tab-blue}{HTML}{1f77b4}
\definecolor{tab-orange}{HTML}{ff7f0e}
\definecolor{tab-green}{HTML}{2ca02c}
\definecolor{tab-red}{HTML}{d62728}
\definecolor{tab-purple}{HTML}{9467bd}
\definecolor{tab-brown}{HTML}{8c564b}
\definecolor{tab-pink}{HTML}{e377c2}
\definecolor{tab-gray}{HTML}{7f7f7f}
\definecolor{tab-olive}{HTML}{bcbd22}
\definecolor{tab-cyan}{HTML}{17becf}

\begin{tikzpicture}[x=70pt,y=70pt]
  \begin{scope}[shift={(20pt, -50pt)}]
  \coordinate (q2a) at (0,5);
  \coordinate (q2b) at (1,5);
  \coordinate (q1a) at (0,4);
  \coordinate (q1b) at (1,4);
  \coordinate (q1c) at (0.5,3);
  \coordinate (q1d) at (0.25,3);
  \coordinate (q1e) at (0.5,2);
  \coordinate (q1f) at (0.75,3);

  \def\points{(q2a), (q2b), (q1a), (q1b), (q1d), (q1f)}
  \def\cpoints{(q1c), (q1e)}

  \draw[tab-gray!50,thin] ([xshift=-40pt]q2a) node[black,anchor=east] {\qtwoqone}                     -- ([xshift=40pt]q2b) node[black, anchor=west] {level $\ell=0$};
  \draw[tab-gray!50,thin] ([xshift=-40pt]q1a) node[black,anchor=east] {$h$ \isoqtwoqone}                        -- ([xshift=40pt]q1b) node[black, anchor=west] {level $\ell=1$};
  \draw[tab-gray!50,thin] ([xshift=-40pt]q1a|-q1d) node[black,anchor=east] {$2h$ \isoqtwoqone} -- ([xshift=40pt]q1b|-q1d) node[black, anchor=west] {level $\ell=2$};
  \draw[tab-gray!50,thin] ([xshift=-40pt]q1a|-q1e) node[black,anchor=east] {$4h$ \isoqtwoqone} -- ([xshift=40pt]q1b|-q1e) node[black, anchor=west] {level $\ell=3$};

  \def\h{2pt};
  \draw[ultra thick, shorten >=\h, ->, >=latex] (q2a) -- (q1a);
  \draw[ultra thick, shorten >=\h, ->, >=latex] (q2b) -- (q1b);

  \draw[tab-blue, ultra thick, shorten >=\h, ->, >=latex] (q1a) -- (q1c);
  \draw[tab-blue, ultra thick, shorten >=\h, ->, >=latex] (q1c) -- (q1b);

  \draw[tab-green, ultra thick, shorten >=\h, ->, >=latex] (q1a) -- (q1d);
  \draw[tab-green, ultra thick, shorten >=\h, ->, >=latex] (q1d) -- (q1e);
  \draw[tab-green, ultra thick, shorten >=\h, ->, >=latex] (q1e) -- (q1f);
  \draw[tab-green, ultra thick, shorten >=\h, ->, >=latex] (q1f) -- (q1b);

  \foreach \p in \points {
    \draw[draw=black!50,fill=tab-gray] \p circle (2pt);
  }
  \foreach \p in \cpoints {
    \draw[draw=black!50,fill=tab-red] \p circle (2pt);
  }

  \draw[|-|,tab-blue] ([xshift=100pt]q1b) --node[midway,fill=white,rotate=90] {two-level} ([xshift=100pt]q1b|-q1d);
  \draw[|-|,tab-green] ([xshift=110pt]q1b) --node[midway,fill=white,rotate=90] {three-level} ([xshift=110pt]q1b|-q1e);

  \draw[thick, tab-gray, shorten >=\h, ->, >=latex] (q2a) -- node[midway,fill=white] {$\gamma=0$} (q2b);
  \draw[thick, tab-gray, shorten >=\h, ->, >=latex] (q1a) -- node[midway,fill=white] {$G^{\gamma}\equiv 0$} (q1b);

  \end{scope}
  \begin{scope}[shift={(240pt, -50pt)}]
  \coordinate (q2a) at (0,5);
  \coordinate (q2b) at (2,5);
  \coordinate (q1a) at (0,4);
  \coordinate (q1b) at (2,4);
  \coordinate (q1a1) at (1,4);
  \coordinate (q1c) at (0.5,3);
  \coordinate (q1c1) at (1.5,3);

  \def\points{(q2a), (q2b), (q1a), (q1b), (q1a1)}
  \def\cpoints{(q1c), (q1c1)}

  \def\h{2pt};
  \draw[ultra thick, shorten >=\h, ->, >=latex] (q2a) -- (q1a);
  \draw[ultra thick, shorten >=\h, ->, >=latex] (q2b) -- (q1b);

  \draw[tab-blue, ultra thick, shorten >=\h, ->, >=latex] (q1a) -- (q1c);
  \draw[tab-blue, ultra thick, shorten >=\h, ->, >=latex] (q1c) -- (q1a1) node[anchor=south,tab-gray] {$\gamma=2$};
  \draw[tab-blue, ultra thick, shorten >=\h, ->, >=latex] (q1a1) -- (q1c1);
  \draw[tab-blue, ultra thick, shorten >=\h, ->, >=latex] (q1c1) -- (q1b);

  \foreach \p in \points {
    \draw[draw=black!50,fill=tab-gray] \p circle (2pt);
  }
  \foreach \p in \cpoints {
    \draw[draw=black!50,fill=tab-red] \p circle (2pt);
  }
  \end{scope}
    \end{tikzpicture}

    \caption{Multilevel scheme for the \qtwoqone{} operator.  Coarse levels are refinements in $h$ on the \isoqtwoqone{} discretization, resulting in either a two-level scheme (blue) or a multilevel scheme (green). On the left, a single V-cycle is executed on the \isoqtwoqone{} problem (with $\gamma=1$); on the right, two V-cycles are executed (with $\gamma=2$). The exact, coarsest level solves are marked in red.}\label{fig:multilevel}
\end{figure}

\subsection{Relaxation}

The relaxation method $M_{\ell}^{-1}$ applied to the saddle-point system $K_{\ell}$ is either chosen as
an additive Vanka~\cite{farrell2019c, PFarrelletal2019a} or Braess-Sarazin~\cite{Braess} coupled relaxation scheme.
In this section, we drop the operator subscript $\ell$, since the relaxation operator construction is the same on all levels of the multigrid hierarchy.

\begin{figure}[!ht]
   \centering
   \newcommand{\mygrid}[3]{      \def\n{#1}
      \def\m{#2}

      \foreach \y in {0,\m,...,\n}{         \draw[gray, loosely dotted] (-1, \y) -- (\n+1, \y);
      }
      \foreach \x in {0,\m,...,\n}{         \draw[black, loosely dotted] (\x, -1) -- (\x, \n+1);
      }

      \foreach \y in {0,2,...,\n}{         \draw[black, thick] (0, \y) -- (\n, \y);
      }
      \foreach \x in {0,2,...,\n}{         \draw[black, thick] (\x, 0) -- (\x, \n);
      }

      \ifnum#3=1,{
         \foreach \y in {1,2,...,\n}{            \draw[black, densely dotted] (0, \y) -- (\n, \y);
         }
         \foreach \x in {1,2,...,\n}{            \draw[black, densely dotted] (\x, 0) -- (\x, \n);
         }
      }
      \else{}\fi

      \def\h{0.15}
      \draw[thick, draw=red!50!red,fill=white] (\n/2-\h, \n/2-\h) rectangle +(2*\h,2*\h);

      \foreach \y in {0,1,...,\n}{         \foreach \x in {0,1,...,\n}{            \pgfmathtruncatemacro\resultX{Mod(\x,2)==0?0:1}
            \pgfmathtruncatemacro\resultY{Mod(\y,2)==0?0:1}
            \pgfmathtruncatemacro\resultXY{\resultX*\resultY}
            \ifnum\resultX=0,{               \ifnum\resultY=0,{\draw[thick, draw=black!60, fill=black!60] (\x, \y) circle (2pt);}
               \else{\node[mark size=2pt,color=black!100, fill=white] at (\x,\y) {\pgfuseplotmark{o}};}\fi            }
            \else{               \ifnum\resultXY=1,{\node[mark size=3pt,color=black!60] at (\x,\y) {\pgfuseplotmark{diamond*}};}
               \else{\ifnum\resultX=1,{\node[mark size=3pt,color=black!60] at (\x,\y) {\pgfuseplotmark{triangle*}};} \else {}\fi}
               \fi
            }\fi
         }
      }
   }
   \subfigure[\qtwoqone]{   \begin{tikzpicture}[x=20pt,y=20pt]
      \begin{scope}[shift={(10pt, 10pt)}]
         \mygrid{4}{2}{0}
      \end{scope}
   \end{tikzpicture}
   }
   \hspace{20mm}
   \subfigure[\isoqtwoqone]{   \begin{tikzpicture}[x=20pt,y=20pt]
      \begin{scope}[shift={(10pt, 10pt)}]
         \mygrid{4}{1}{1}
      \end{scope}
   \end{tikzpicture}
   }
   \caption{At left, construction of a typical Vanka patch for the \qtwoqone{} discretization.  At right, construction of a typical Vanka patch for the \isoqtwoqone{} discretization.  Note that the patch in (b) includes one pressure DoF and all the velocity DoFs that geometrically lie within (or on the boundary) of the sub-domain associated with the elements that include that pressure DoF.}\label{fig:Vanka_patches}
\end{figure}
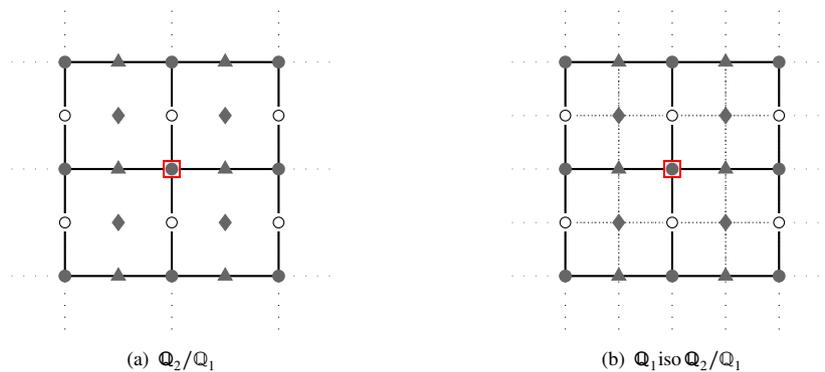

Letting $n_p$ be the number of pressure degrees of freedom (DoFs), Vanka relaxation partitions the $n\times n$ system $K$ into $n_p$ overlapping saddle-point problems, consisting of $2\times2$ patches of elements around each nodal pressure DoF.  An important feature of the construction of Vanka relaxation (in contrast to so-called ``star'' relaxation~\cite{farrell2019pcpatch}) is that each patch consists of all velocity degrees of freedom that are on the closure of the elements adjacent to the nodal pressure DoF, and not just those in the interior of mesh entities (elements and faces) adjacent to the node.  Each patch contains a single pressure DoF, however.  \cref{fig:Vanka_patches} shows the construction of a Vanka patch for both the \qtwoqone{} and \isoqtwoqone{} discretizations; it is important to note that the construction of the Vanka patches in the \isoqtwoqone{} discretization is based on the (coarse) pressure elements and not the (fine) velocity elements.
For each patch, indexed by pressure DoF $i$, we form a (binary) restriction operator, $V_i$, which selects those entries in a global vector that appear on patch $i$.  The system matrix is then projected onto the patch DoFs by a triple matrix product, $V_i K V_i^T$. A single iteration of Vanka is then given by
\begin{equation*}\label{ASM-precondition}
  M_V^{-1} = \sum_{i=1}^{n_p} V_i^{T} W_i \left (V_i K V_i^T \right )^{-1}  V_i,
\end{equation*}
where the diagonal weighting matrix $W_i$ is defined such that each diagonal entry is equal to the reciprocal of the number of patches that
contain the associated degree of freedom.

In contrast to Vanka relaxation, Braess-Sarazin relaxation retains the block structure of the Stokes system to produce a suitable relaxation method. A single application of Braess-Sarazin relaxation is formulated as the (approximate) solution to the following system
\begin{equation}\label{eq:braess}
    M_{BS} \delta x =
    \begin{bmatrix}
        \alpha \tilde{A}  &   B^T    \\
           B  &  0
    \end{bmatrix}
     \begin{bmatrix}
       \delta \vec{u}    \\
       \delta p
    \end{bmatrix}
    =
     \begin{bmatrix}
       r_{\vec{u}}    \\
       r_p
    \end{bmatrix},
\end{equation}
where $\alpha>0$ is a relaxation parameter, $\tilde{A}$ is some approximation to $A$, $r_{\vec{u}}$ and $r_p$ are the components of the current residual, and $\delta \vec{u}$ and $\delta p$ are the components of the correction. The system in~\cref{eq:braess} is solved in two sequential steps
\begin{subequations}\label{eq:braess_componentwise}
\begin{align}
   (B \tilde{A}^{-1}B^T) \delta p &=   B\tilde{A}^{-1}r_{\vec{u}}- \alpha {\hskip.02in} r_p, \label{eq:bs1} \\
  \delta \vec{u}  &= \frac{1}{\alpha} \tilde{A}^{-1} (r_{\vec{u}}-B^T \delta p).  \label{eq:bs2}
\end{align}
\end{subequations}
Different Braess-Sarazin variations
are devised using different choices for $\tilde{A}$ and whether or not~\cref{eq:bs1} is solved \textit{exactly}.
A common approach is to take $\tilde{A}$ to be the diagonal matrix
defined by $\tilde{A}_{ii} = A_{ii}$ for all $i$.
This allows $B\tilde{A}^{-1}B^T$ to be directly computed, as well as matrix-vector products with $\tilde{A}^{-1}$.
Note that the exact solution of \cref{eq:bs1} is impractical and, so, we consider the inexact Braess-Sarazin relaxation (IBSR), where \cref{eq:bs1} is approximately solved via a single weighted Jacobi iteration, with a weighting factor $\beta$.

Both Vanka and IBSR relaxation include several parameters that must be chosen.  While parameter choices are known in some cases~\cite{he2019local, PFarrelletal2019a}, our defect-correction framework is somewhat unique and so we employ
LFA to inform our parameter choices for the relaxation methods and for the MG cycle defined by~\cref{eq:TG-operator,eq:TG-operator-full}.

\section{Local Fourier Analysis}\label{sec:LFA}

LFA is a well-known and valuable tool in predicting and analyzing algorithmic performance for the solution of discretized PDEs~\cite{Trottenberg2001,wienands2004practical}. In particular, LFA is often used to guide parameter choices for
multigrid components,
including those within relaxation schemes and grid-transfer operators. Here, we apply LFA to the monolithic multigrid methods in~\cref{sec:multigrid} and aim to optimize the spectral radius of the two-grid error propagation operator~\cref{eq:TG-operator} over the possible choices of parameters.

\subsection{Definition and notations}\label{sec:Fourier_definitions}

For completeness, we give a brief introduction to LFA~\cite{Trottenberg2001,wienands2004practical}. First, consider a two-dimensional infinite uniform grid
\begin{equation*}
  \mathcal{G}=\big\{\vec{x}:=(x_1,x_2)=(k_{1},k_{2})h,\quad (k_1,k_2)\in \mathbb{Z}^2\big\},
\end{equation*}
with uniform grid size $h$.
Let $A$ be a scalar Toeplitz operator defined by entries in a ``stencil'', $s_{\boldsymbol{\kappa}}\in \mathbb{R}$, where
$\boldsymbol{\kappa} \in V\subset \mathbb{Z}^2$ is a finite index set over which the stencil is nonzero,
that acts on a vector $w(\boldsymbol{x}) \in l^2(\mathcal{G})$ as follows:
\begin{equation*}\label{defi-symbol-classical}
    Aw(\boldsymbol{x})=\sum_{\boldsymbol{\kappa}\in V}s_{\boldsymbol{\kappa}}w(\boldsymbol{x}+\boldsymbol{\kappa}h).
\end{equation*}
Such operators can be formally diagonalized by the Fourier modes $\psi(\boldsymbol{\theta},\boldsymbol{x})= e^{\iota\boldsymbol{\theta}\cdot\boldsymbol{x}/{h}}=e^{\iota \theta_1x_1/h}e^{\iota \theta_2x_2/h}$, where $\boldsymbol{\theta}=(\theta_1,\theta_2)$ and $\iota^2=-1$.  Thus, we use $\psi(\boldsymbol{\theta},\boldsymbol{x})$ as a Fourier basis with $\boldsymbol{\theta}\in \left[-\frac{\pi}{2},\frac{3\pi}{2}\right)^{2}$. Considering standard coarsening by a factor of 2 in each direction, the relevant spaces of low and high frequencies  are given by
\begin{equation*}
  \boldsymbol{\theta}\in T^{{\rm low}} =\left[-\frac{\pi}{2},\frac{\pi}{2}\right)^{2}, \, \boldsymbol{\theta}\in T^{{\rm high}} =\displaystyle \left[-\frac{\pi}{2},\frac{3\pi}{2}\right)^{2} \bigg\backslash \left[-\frac{\pi}{2},\frac{\pi}{2}\right)^{2}.
\end{equation*}

LFA provides two possible predictions of multigrid performance through the so-called ``smoothing'' and ``two-grid convergence'' factors, which often offer
sharp predictions of actual multigrid performance. Unfortunately, the simpler LFA smoothing factor can provide poor predictions when used on complicated or higher-order operators~\cite{HM2020LFALaplace}. Thus, we instead focus here on the two-grid LFA convergence factor, which includes the full details of the coarse-grid correction process.  To do this, we define the following harmonic modes
\begin{equation*}
\boldsymbol{\theta}^{\boldsymbol{\xi}}=(\theta_1^{\xi_1},\theta_2^{\xi_2})=\boldsymbol{\theta}+\pi\cdot\boldsymbol{\xi},\,\,
\boldsymbol{\theta}:=\boldsymbol{\theta}^{00}\in T^{{\rm low}},
\end{equation*}
where $\boldsymbol{\xi}=(\xi_1,\xi_2)\in\big\{(0,0),(1,0),(0,1),(1,1)\big\}$.

\begin{definition}\label{formulation-symbol-classical}
For a given Toeplitz operator, $A$, $\widetilde{A}(\boldsymbol{\theta})=\displaystyle\sum_{\boldsymbol{\kappa}\in V}s_{\boldsymbol{\kappa}}e^{\iota \boldsymbol{\theta}\cdot\boldsymbol{\kappa}}$ is the symbol of $A$.
\end{definition}
Note that for all Fourier modes, $\psi(\boldsymbol{\theta},\boldsymbol{x})$, we have
\[
A\psi(\boldsymbol{\theta},\boldsymbol{x})= \widetilde{A} (\boldsymbol{\theta})\psi(\boldsymbol{\theta},\boldsymbol{x}).
\]
If the multigrid relaxation scheme, $M$, is also represented by a Toeplitz operator (or its inverse), then the symbol for the error-propagation operator $\widetilde{S}(\boldsymbol{\theta}) = I-\widetilde{M}^{-1} (\boldsymbol{\theta})\widetilde{A} (\boldsymbol{\theta})$ provides information on how the relaxation scheme damps errors at each Fourier frequency.  The LFA smoothing factor arises from computing the maximum (absolute) value of $\widetilde{S}(\boldsymbol{\theta})$ over the high-frequency set, $T^{{\rm high}}$.

To get a better picture of multigrid convergence, we study how the relaxation scheme represented by $\widetilde{S}(\boldsymbol{\theta})$ interacts with the coarse-grid correction.  Here, we must account for the coupling of fine-grid errors in the coarse-grid correction process.
For each low-frequency mode, $\boldsymbol{\theta}\in T^{\rm low}$, we define a four-dimensional harmonic space,
\begin{equation*}\label{eq:simple-invariant-space}
  \vec{F}(\boldsymbol \theta)={\rm span}\Big\{ \psi(\boldsymbol{\theta^{\xi}},\cdot): \boldsymbol{\xi}\in\big\{(0,0),(1,0),(0,1),(1,1)\big\}\Big\},
\end{equation*}
which is invariant for standard full-coarsening two-grid algorithms under certain assumptions on their grid-transfer operators. The \textit{symbol}  of the two-grid algorithm
is defined analogously to~\cref{formulation-symbol-classical}, accounting for how interpolation and restriction map between grids.  For the given interpolation and restriction operators, $P$ and $R$, we can define their symbols by how they map harmonic modes between grids.  For interpolation, a single coarse-grid mode with frequency $2\boldsymbol{\theta}$ is mapped onto $\vec{F}(\boldsymbol \theta)$, resulting in a $4\times 1$ symbol for $P$, denoted $\widetilde{P}_{2g}$, with entries corresponding to each harmonic frequency.  The restriction is similarly mapped into a $1\times 4$ symbol, $\widetilde{R}_{2g}$, representing how each fine-grid harmonic frequency is mapped onto the coarse-grid mode with frequency $2\boldsymbol{\theta}$.  As a result, the symbol of the two-grid error propagation operator is a $4\times 4$ matrix given by
\begin{equation}\label{eq:LFA-represenation}
 \widetilde{E}_{2g}(\boldsymbol{\theta})= \widetilde{S}_{2g}^{\nu_2}(\boldsymbol{\theta})
 \left(I-\widetilde{{P}}_{2g}(\boldsymbol{\theta})(\widetilde{{A}}_{C}(2\boldsymbol{\theta}))^{-1}
         \widetilde{R}_{2g}(\boldsymbol{\theta})
         \widetilde{A}_{2g}(\boldsymbol{\theta})
\right)
\widetilde{S}_{2g}^{\nu_1}(\boldsymbol{\theta}),
\end{equation}
where
\begin{align*}
\widetilde{A}_{2g}(\boldsymbol{\theta})&=\text{diag}\left\{\widetilde{{A}}(\boldsymbol{\theta}^{00}), \widetilde{{A}}(\boldsymbol{\theta}^{10}),\widetilde{{A}}(\boldsymbol{\theta}^{01}),
\widetilde{{A}}(\boldsymbol{\theta}^{11})\right\}, \\
\widetilde{{R}}_{2g}(\boldsymbol{\theta})&= \left[\widetilde{R}(\boldsymbol{\theta}^{00}),\widetilde{R}(\boldsymbol{\theta}^{10}),
\widetilde{R}(\boldsymbol{\theta}^{01}),\widetilde{R}(\boldsymbol{\theta}^{11})\right], \\
\widetilde{P}_{2g}(\boldsymbol{\theta})&=\left[\widetilde{P}(\boldsymbol{\theta}^{00}), \widetilde{P}(\boldsymbol{\theta}^{10}),
\widetilde{P}(\boldsymbol{\theta}^{01}), \widetilde{P}(\boldsymbol{\theta}^{11}) \right]^T, \\
\widetilde{S}_{2g}(\boldsymbol{\theta})&=\text{diag}\left\{\widetilde{S}(\boldsymbol{\theta}^{00}),
\widetilde{S}(\boldsymbol{\theta}^{10}),\widetilde{S}(\boldsymbol{\theta}^{01}),
\widetilde{S}(\boldsymbol{\theta}^{11})\right\}.
\end{align*}
Here, ${\rm diag}\{A_1,A_2,A_3,A_4\}$ is the diagonal matrix with diagonal entries, $A_1, A_2, A_3$, and $A_4$, and $\widetilde{A}_C$ is the symbol of the coarse-grid operator.
\begin{definition}
The two-grid LFA convergence factor, $\hat{\rho}$, is defined as
\begin{equation}\label{eq:rho-definition}
  \hat{\rho} = \sup_{\boldsymbol{\theta}\in T^{\text{low}}} \rho\left(\widetilde{E}_{2g}(\boldsymbol{\theta})\right),
\end{equation}
where $\rho( \widetilde{E}_{2g}(\boldsymbol{\theta}))$ denotes the spectral radius of matrix $ \widetilde{E}_{2g}(\boldsymbol{\theta})$.
\end{definition}
In this work,  we will apply LFA to the monolithic multigrid algorithms described above, and optimize the two-grid convergence factor  using a robust optimization framework recently developed for cases such as this with multiple parameters to be optimized~\cite{LFAoptAlg2020}.  While we do not compute $\hat{\rho}$ from~\Cref{eq:rho-definition} exactly, we use this to denote the values obtained from the optimization algorithm that, necessarily, samples the spectral radius of $\widetilde{E}_{2g}(\boldsymbol{\theta})$ at only a finite number of points.

\subsection{Fourier representation for
\texorpdfstring{\qtwoqone}{Q2/Q1} and
\texorpdfstring{\isoqtwoqone}{Q1isoQ2/Q1} discretizations}

LFA for coupled systems discretized using mixed finite-element approaches is necessarily more complicated than in the scalar case.  While the staggered finite-difference case can be successfully handled using slight generalizations of the scalar case~\cite{ANiestegge_KWitsch_1990a, Siv_1991a, Trottenberg2001}, greater adaptation is required for higher-order or otherwise non-nodal finite-element discretizations~\cite{SPMacLachlan_CWOosterlee_2011a, CRodrigo_etal_2016a, he2019local, HM2020LFALaplace}.  The key step is in realizing that the resulting discrete operators (including the infinite-grid analog of~\cref{equation:saddle}) can be reordered into block structured linear systems with Toeplitz blocks.  Once suitably reordered, each block can be diagonalized using the classical Fourier approach.  In this way, LFA of more complicated operators in a two-grid method results in block Fourier symbols that can be assembled as in~\cref{eq:LFA-represenation} and maximized to yield a two-grid LFA convergence factor, as in~\Cref{eq:rho-definition}.

For the \qtwoqone{} and \isoqtwoqone{} systems considered here, the natural structure is as $9\times 9$ block systems.  Consider the velocity DoFs as pictured in either Figure~\ref{fig:meshes} or Figure~\ref{fig:Vanka_patches}; both the \qtwo{} and \isoqtwo{} spaces have basis functions associated with the nodes of the quadrilateral mesh, as well as along horizontal and vertical mesh edges, and at the center of the mesh cells.  Since there are 2 components (horizontal and vertical) to the velocity, this gives eight {\it types} of DoFs associated with the velocity discretization.  The $\mathbb{Q}_1$ pressure gives the ninth set of DoFs in the block $9\times 9$ structure.  A key realization necessary for the LFA is that, when the system matrices (or their infinite-grid analogs) are reordered blockwise, each block in the reordered system has Toeplitz structure.  Thus, we can diagonalize each block in the system by a Fourier transform, as discussed in~\Cref{sec:Fourier_definitions}.  After diagonalizing the blocks, the system can be reordered frequency-wise, into a block operator with $9\times 9$ blocks, by collecting all rows/columns corresponding to each Fourier frequency into adjacent rows/columns of the operator.  If we further order the frequencies in harmonic sets, we get a block $36\times 36$ structure associated with each harmonic space, $\vec{F}(\boldsymbol \theta)$.  While analysis of the relaxation operators is most natural using the $9\times 9$ block structure, complete analysis of the two-grid convergence factor requires use of the $36\times 36$ block structure.

The block symbol for the \qtwoqone{} Stokes discretization was derived in earlier work~\cite{he2019local}, as are symbols for the \qtwo{} Laplacian and the canonical finite-element interpolation and restriction operators for the \qtwo{} finite-element space~\cite{HM2020LFALaplace}.  Symbols for both Braess-Sarazin~\cite{he2019local} and additive Vanka~\cite{PFarrelletal2019a} relaxation for this discretization were also computed in earlier work.  There is, in principle, very little different to computing these symbols for the \isoqtwoqone{} discretization.
The reordered block structure of the  \isoqtwoqone{} saddle point matrix is the same as that of the \qtwoqone{} matrix. Once again, each individual block of the $9 \times 9 $ system is a Toeplitz matrix, though the specific entries are now different than for the \qtwoqone{} case. For example, a block row associated with one {\it type} of velocity DoF would have five nonzero blocks. If this {\it type} corresponds to horizontal velocity DoFs located at pressure cell centers (denoted by a diamond marker in Figure~\ref{fig:meshes}), then four of these nonzero blocks correspond to {\it types} associated with horizontal velocity DoFs and the fifth block corresponds to the gradient operator acting on the pressure.  The diagonal block within this block row is itself diagonal, reflecting the fact that there are no matrix connections between different velocity {\it types} in the \isoqtwo{} Laplacian matrix. Referring to the right side of Figure~\ref{fig:meshes}, this corresponds to the fact that two adjacent velocity DoFs in the velocity mesh
always have different markers.  Two of the other nonzero velocity blocks within this same block row are Toeplitz matrices with only 2 nonzeros per row.
Again referring to the right side of Figure~\ref{fig:meshes}, the nonzeros in one of these blocks correspond to stencil entries associated with the North
and South horizontal velocity DoFs, located on the horizontal edges of the mesh, which share the same marker. The two nonzeros in the other block correspond to stencil entries associated with the East
and West horizontal velocity DoFs, now on vertical edges of the mesh, which again share the same marker. The final nonzero velocity block in this same block row would have four nonzeros that correspond
to stencil entries associated with the corner horizontal velocity DoFs, located at the nodes of the mesh, which all share the same marker. The Toeplitz block for the gradient operator mapping onto these velocity DoFs would have four nonzero entries
per row, associated with the four nodal pressure DoFs for the cell.
It should be noted that the number of nonzeros per row in a gradient Toeplitz block depends on the specific block row that is being considered, i.e. the
particular velocity marker associated with the block row in Figure~\ref{fig:meshes}.
For the stencil representation of each Toeplitz matrix,
standard Fourier techniques, as described above, give the block-structured symbol of the operator.  Similarly, the symbols for interpolation and restriction can be computed~\cite{HM2020LFALaplace} in the \isoqtwo{} setting, and those for Braess-Sarazin~\cite{he2019local} and additive Vanka~\cite{PFarrelletal2019a} relaxation as well.
Since the full calculations are rather tedious, we omit the details, noting that there are existing software packages~\cite{MBolten_HRittich_2018a, KKahl_NKintscher_2020a} that could be used to compute the symbols automatically.

\section{Numerical Results}\label{sec:num}

To evaluate the effectiveness of \isoqtwoqone{} multigrid preconditioners for the \qtwoqone{} system,
we start by comparing LFA predicted two-grid convergence factors against measured convergence factors
for the Stokes system with Dirichlet boundary conditions on a uniform square mesh with $h=\frac{1}{64}$.
For convenience, the \qtwoqone{} discretization matrices are assembled using Firedrake~\cite{FRathgeber_etal_2017a, kirby2018solver}, while the \isoqtwoqone{} discretization matrices are assembled by forming a finer-mesh \qtwoqone{} discretization, then applying coarsening-in-order to the velocity DoFs and coarsening-in-space to the pressure DoFs.  The defect-correction preconditioner is implemented as a GMG variant taken from the PyAMG library~\cite{OlSc2018}, with our own implementation of the relaxation schemes.
The measured convergence factors are then computed by running the solver $100$ times using a zero right-hand side and a different random initial guess each time. The vector $\ell^2$ norm of the last $\min(m_s,10)$ residuals are collected for each of the $100$ runs, where $m_s$ is the
number of iterations required to reduce the residual by a factor of $10^{-10}$ for the $s^{th}$ solve.  When this criteria
is not met in $100$ iterations, we take $m_s = 100$.  Collected residual norms are then fit by a linear least-squares model of the form $\xi_0 + \xi_1 j = \ln r_{js}$ where $r_{js}$ is the $j^{th}$ residual within the $\min(m_s,10)$ collected residuals from the $s^{th}$ run. The averaged asymptotic convergence factor is then taken as $\rho = e^{\xi_1}$.

Within the multigrid cycle, the number of pre- and post-relaxation sweeps is always $1$ for $\ell \ge 1$. To distinguish between different preconditioners based on~\cref{eq:TG-operator-full}, we focus on the triple $(\nu_1, \nu_2, \gamma)$, where $\nu_{1}$ and $\nu_2$ are the number of pre-/post- relaxation sweeps on the $K_0$ system, and $\gamma$ is the number of $h$-multigrid cycles on the $K_1$ system.  For fixed values of these parameters, we use the LFA optimization framework of Brown et al.~\cite{LFAoptAlg2020} to choose the remaining parameters in the method.  Referring to~\cref{eq:TG-operator-full}, these parameters include a damping parameter for the \isoqtwoqone{} cycling, $\tau$, outer relaxation damping parameters for both the \qtwoqone{} and \isoqtwoqone{} relaxation schemes and, if IBSR relaxation is used, inner relaxation parameters for that relaxation.  The notation for these parameters is summarized in~\cref{table:param-ref}.
\begin{table}[!ht]
\centering
\begin{tabular}{cl}
\toprule
   Symbol &Description\\
\midrule
  $K_0$         & saddle point system for \qtwoqone\\
  $K_1$         & saddle point system for \isoqtwoqone\\
  $K_\ell$      & coarse saddle point \isoqtwoqone{} system for $\ell \ge 2$ \\
  $M_\ell$      & relaxation operator on level $\ell$ \\
\midrule
  $\nu_1$       & number of pre-relaxation sweeps used on the $K_0$ system \\
  $\nu_2$       & number of post-relaxation sweeps used on the $K_0$ system \\
  $\gamma$      & number of $h$-multigrid cycles used to solve the $K_1$ system \\
\midrule
  $\tau$        & \isoqtwoqone-cycle solution damping parameter on level $1$ \\
  $\alpha_0$    & Inexact Braess-Sarazin (IBS) relaxation parameter on level $0$ \\
  $\beta_0$     & weighted Jacobi relaxation parameter in IBSR  on level $0$ \\
  $\omega_0$    & relaxation global update parameter on level $0$ \\
  $\alpha_1$    & Braess-Sarazin (BS) relaxation parameter for all $\ell \ge 1$\\
  $\beta_1$     & BS weighted Jacobi relaxation parameter  for all $\ell \ge 1$\\
  $\omega_1$    & relaxation global update parameter   for all $\ell \ge 1$\\
\midrule
  $\hat{\rho}$  & LFA predicted two-grid convergence factor \\
  $\rho$        & measured convergence factor for Dirichlet BC problem \\
  $m$           & measured iteration count for multigrid or preconditioned FGMRES\\
\midrule
  superscript P & problem with periodic boundary conditions \\
\bottomrule
\end{tabular}
  \caption{Notation Reference Table}\label{table:param-ref}
\end{table}

\subsection{Two-grid Results}

We start by verifying that the measured two-grid convergence factors, $\rho$, agree with the LFA predictions, $\hat{\rho}$. \Cref{table:LFA-mixed-two-grid1} demonstrates that there is generally good agreement. As expected, we also see that increasing the number of $K_0$ relaxation sweeps improves convergence. For example, increasing the number of IBSR sweeps from $(\nu_1,\nu_2, \gamma)=(1,0,1)$ to $(\nu_1,\nu_2,\gamma)=(2,2,1)$ improves the measured convergence factor from \highlight{blue!20}{blue!20}{0.16} to \highlight{blue!20}{blue!20}{0.05}. Similar V-cycles for Vanka relaxation demonstrate an analogous trend, where the convergence factor decreases from \highlight{blue!20}{blue!20}{0.31} to \highlight{blue!20}{blue!20}{0.15}.
Interestingly, there is some notable variation in both the damping parameter, $\tau$, and the outer relaxation weights $\omega_\ell$.  In results not reported here,
the optimal value for $\tau$ was found to be $1/2$ for a simpler preconditioner with error-propagation operator $I-\tau K_1^{-1}K_0$ (i.e., taking $\nu_1 = \nu_2 = 0$ and $\gamma \rightarrow \infty$).  Here, we see that adding relaxation on the $K_0$ system and using a more practical approximation to $K_1^{-1}$ result in values for $\tau$ much closer to 1.  Depending on the $(\nu_1,\nu_2,\gamma)$ values, we see values for $\omega_\ell$ varying in the range from roughly $1/2$ to $1$, and similar variations in the inner IBSR parameters.
\begin{table}[!ht]
  \caption{LFA predictions ($\hat{\rho}$) compared with the measured convergence factors ($\rho$) for
inexact Braess-Sarazin (BS) and Vanka (V) relaxation.  As $\alpha_\ell$ and $\beta_\ell$ are only needed for IBSR relaxation, we do not record values of these parameters for Vanka relaxation.}
\centering
\begin{tabular}{l ccc |rrrrrrr|rr}
\toprule
   &$\nu_1$&$\nu_2$&$\gamma$ &$\tau$ &$\omega_0$ &$\omega_1$ &$\alpha_0$ &$\beta_0$ &$\alpha_1$ &$\beta_1$ &$\hat{\rho}$ &$ \rho$\\
\midrule
\multirow{4}{5mm}{BS} &1&0&1  &.87   &1.02  &.90   &1.20  &.79  &.74 &.88 &.18  &\highlight{blue!20}{blue!20}{.16}\\
                      &1&0&2  &.86   &.75   &1.04  &.99   &.70  &1.15&.88 &.20  &.18\ \\
                      &1&1&1  &.97   &.91   &1.04  &1.02  &.93  &.78 &.86 &\highlight{white}{cyan!100}{.09}  &.11\\
                      &2&2&1  &1.04  &.55   &.59   &.72   &.99  &.85 &1.30&.04  &\highlight{blue!20}{blue!20}{.05}\\
\midrule
\multirow{6}{5mm}{V}  &1&0&1  &.86   &.78   &1.01  &-     &-    &-   &-   &.37  &\highlight{blue!20}{blue!20}{.31}\\
                      &1&0&2  &.90   &.98   &.83   &-     &-    &-   &-   &.21  &.20\\
                      &1&0&3  &.85   &1.02  &.82   &-     &-    &-   &-   &.20  &.19\\
                      &1&1&1  &.99   &.67   &.93   &-     &-    &-   &-   &.29  &\highlight{white}{green!100}{.27}\\
                      &1&1&2  &1.05  &.74   &.78   &-     &-    &-   &-   &.12  &.09\\
                      &2&2&1  &0.98  &.71   &1.05   &-     &-    &-   &-   &.16  &\highlight{blue!20}{blue!20}{.15}\\
\bottomrule
\end{tabular}\label{table:LFA-mixed-two-grid1}
\end{table}

In order to test whether relaxation on both the $K_0$ and $K_1$ systems is necessary for robust convergence, we eliminate the relaxation on one of these systems by setting either $\omega_0=0$ or $\omega_1=0$.
The results in~\cref{table:LFA-mixed-two-grid2} suggest that removing the $K_0$ relaxation is more detrimental to convergence than removing relaxation on the $K_1$ system.  Specifically, we compare $(\nu_1,\nu_2,\gamma) = (0,0,1)$ and $\omega_0 = 0$ with $(\nu_1,\nu_2,\gamma) = (1,1,1)$ and $\omega_1 = 0$.
In this case, there is a single pre- and post-relaxation in the cycle, on either the \isoqtwoqone{} or \qtwoqone{} system.  Similarly, we compare $(\nu_1,\nu_2,\gamma) = (0,0,2)$ and $\omega_0 = 0$ with $(\nu_1,\nu_2,\gamma) = (2,2,1)$ and $\omega_1 = 0$. Here, there are now two pre- and post-relaxation sweeps in the cycle on either the \isoqtwoqone{} or \qtwoqone{} system,
although using $\gamma=2$ results in one additional coarse-grid correction in comparison to the $\gamma=1$ case.  We emphasize that setting $\omega_1 = 0$ omits only relaxation on the ``fine grid'' $K_1$ system in
the two-level method.
In the case of Braess-Sarazin with $\omega_0 = 0$ and $\gamma=2$, the predicted convergence factor is \highlight{green!20}{green!20}{0.52} compared with \highlight{green!20}{green!20}{0.11} for relaxation only on $K_0$ using $(\nu_1,\nu_2,\gamma)=(2,2,1)$. Similar V-cycles for Vanka relaxation demonstrate less of a convergence factor reduction from \highlight{green!20}{green!20}{0.57} to \highlight{green!20}{green!20}{0.29}.   We also note that some of the measured convergence factors \highlight{white}{red}{deviate} from the LFA predications in the case of Braess-Sarazin, while measured Vanka results agree  closely with the LFA predictions.  The two-grid convergence factors for periodic boundary condition problems are expected to match the predicted values exactly~\cite{Trottenberg2001, stevenson1990validity, CRodrigo_etal_2019a}, while sometimes a gap between the LFA predicted factors and the measured factors is observed for the Dirichlet boundary condition case~\cite{he2019local, Trottenberg2001}.
\begin{table}[!ht]
  \caption{Effect of omitting relaxation on $K_0$ and $K_1$. $\omega_0=0$ corresponds only relaxation on the \isoqtwoqone{} system, while $\omega_1=0$ corresponds to only relaxation on the \qtwoqone{} system.}
\centering
\begin{tabular}{l ccc |rrrrrrr|rr}
\toprule
                 &$\nu_1$&$\nu_2$&$\gamma$ &$\tau$ &$\omega_0$ &$\omega_1$ &$\alpha_0$ &$\beta_0$ &$\alpha_1$ &$\beta_1$ &$\hat{\rho}$  &$\rho$\\
 \midrule
\multirow{4}{5mm}{BS}  &0&0&1   &.93   &0       &.24   &-    &-     &.67     &1.56   &.60      &.65 \\
                       &0&0&2   &.63   &0       &.40   &-    &-     &.52     &1.35   &\highlight{green!20}{green!20}{.52}      &.51 \\
                       &1&1&1   &.77   &.45     &0     &.46  &1.37  &-       &-      &.32      &\highlight{white}{red}{.52}\\
                       &2&2&1   &.99   &.85     &0     &.90  &1.12  &-       &-      &\highlight{green!20}{cyan!100}{.11}      &\highlight{white}{red}{.41}\\
\midrule
\multirow{4}{5mm}{V}   &0&0&1   &.64   &0       &.98   &-    &-     &-       &-      &.63      &.65 \\
                       &0&0&2   &.66   &0       &.63   &-    &-     &-       &-      &\highlight{green!20}{green!20}{.57}      &.55 \\
                       &1&1&1   &.48   &.81     &0     &-    &-     &-       &-      &.56      &.51 \\
                       &2&2&1   &.80   &.85     &0     &-    &-     &-       &-      &\highlight{green!20}{green!20}{.29}      &\highlight{white}{green!100}{.25} \\
                       \bottomrule
\end{tabular}\label{table:LFA-mixed-two-grid2}
\end{table}

Considering the information in~\cref{table:LFA-mixed-two-grid1,table:LFA-mixed-two-grid2}, we compare whether relaxation on both $K_0$ and $K_1$ is necessary or  relaxation on $K_0$ system is sufficient. For example, Vanka relaxation with $(\nu_1,\nu_2,\gamma)=(1,1,1)$ in~\cref{table:LFA-mixed-two-grid1} relaxes on both systems, while taking $(\nu_1,\nu_2,\gamma)=(2,2,1)$ with $\omega_1 = 0$ in~\cref{table:LFA-mixed-two-grid2} use the same total number of relaxation sweeps, but only on the $K_0$ system. Their convergence factors are comparable, giving \highlight{white}{green!100}{0.27} and \highlight{white}{green!100}{0.25}, respectively. For Braess-Sarazin, the LFA predicted performance for the same cycle types is also comparable, giving \highlight{white}{cyan!100}{0.09} and \highlight{white}{cyan!100}{0.11}, respectively. However, the measured Braess-Sarazin convergence factors deviate quite significantly from these predictions.  Since there appears to be little benefit to omitting this relaxation even in the ideal setting of periodic boundary conditions, and potentially significant degradation in performance when used with Dirichlet boundary conditions, we do not consider the methods from~\cref{table:LFA-mixed-two-grid2} further and, instead, focus on methods where we allow relaxation on both $K_0$ and $K_1$.

\begin{remark}[Comparison with direct application of GMG to the \qtwoqone{} discretization]
 The $hp$-multigrid convergence rates that we report are quite similar to the  \qtwoqone{} $h$-multigrid results reported in earlier work~\cite{farrell2020local, he2019local}.
 For Vanka relaxation, our two-grid cycle $(1,1,1)$ has the same total number of fine-level relaxation sweeps (considering both $K_0$ and $K_1$) as a $(2,2)$ $h$-multigrid cycle, where a total of two pre- and two post- relaxation sweeps are performed. The measured convergence factor for a unit square domain with periodic boundary conditions is reported as 0.29 for the $h$-multigrid solver~\cite{farrell2020local},
which exactly matches our results.
 For IBSR based GMG, a measured convergence factor of 0.09 is reported~\cite{he2019local} with an $h$-multigrid $(2,2)$ cycle for a problem on a unit square domain with Dirichlet boundary
conditions, which closely matches our rate of .11 for a $(1,1,1)$ cycle in~\cref{table:LFA-mixed-two-grid2}. Here, small differences could potentially arise due to the different number of weighted-Jacobi sweeps in the computation of the pressure correction (two sweeps for the $h$-multigrid benchmark versus one for our method).  Overall, however it is clear that comparable convergence rates are obtained when either GMG is applied directly to the \qtwoqone{} discretization or when it is applied indirectly via an intermediate \isoqtwoqone{} system.
 \end{remark}

\begin{figure}[!ht]
    \centering
    \includegraphics[width=0.4\textwidth]{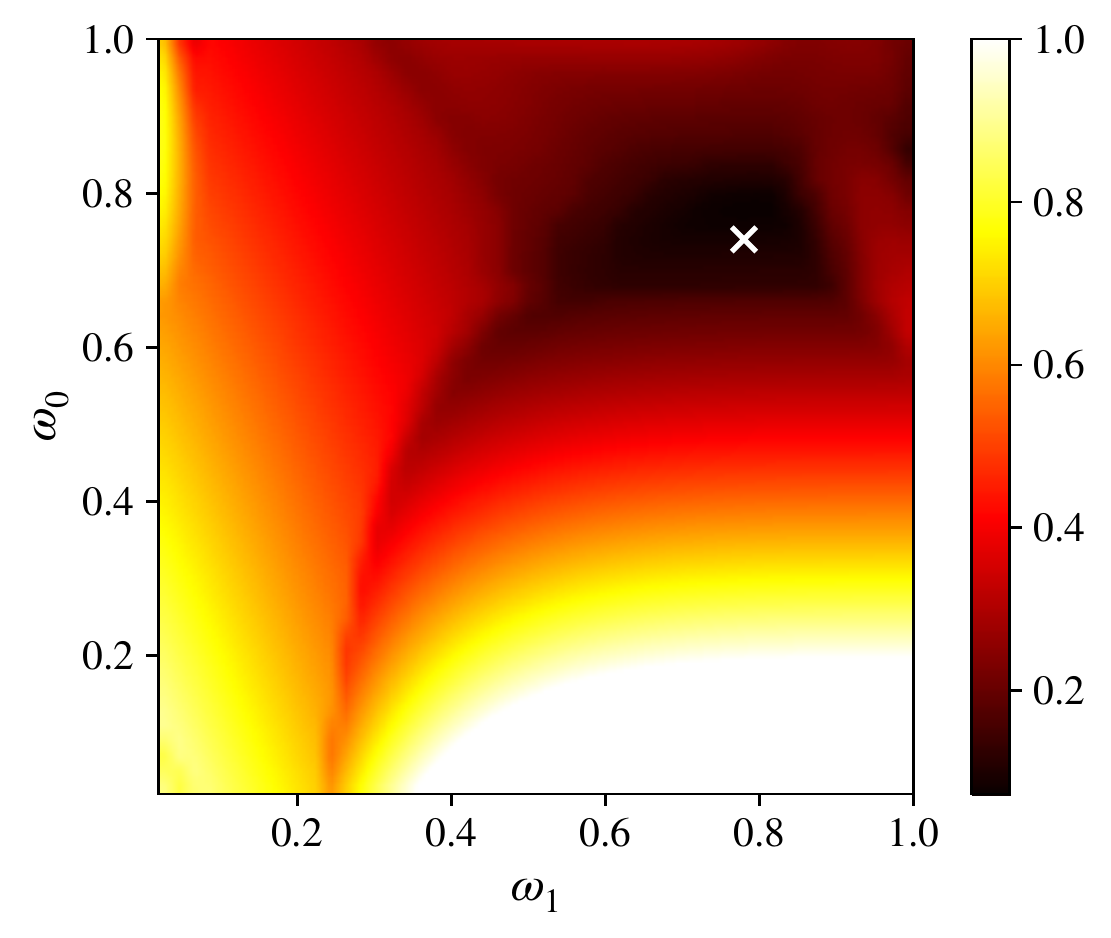}
    \includegraphics[width=0.4\textwidth]{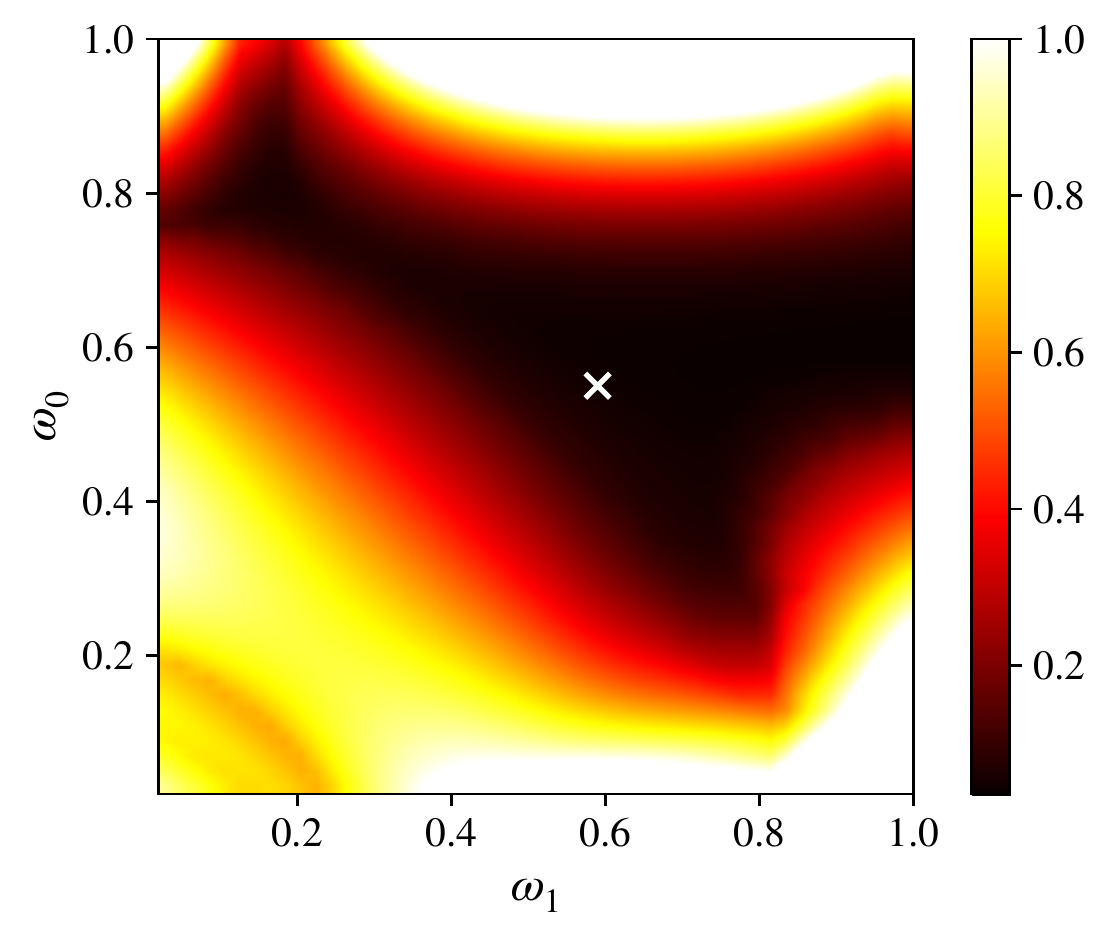}
    \caption{Sensitivity analysis of the measured  two-grid convergence factors for V$(1,1,2)$  (left) and BS$(2,2,1)$ (right) cycles. }\label{fig:sensitivity-plots}
\end{figure}

In order to study the convergence rate sensitivity to parameter choices, we consider the two-grid convergence factor for the same Dirichlet boundary value problem as in the previous section. In case of Vanka relaxation, we consider the V$(1,1,2)$ cycle from~\cref{table:LFA-mixed-two-grid1}, where we fix $\tau$ and vary $\omega_0$ and $\omega_1$ from 0.02 to 1.0. In the left plot of~\cref{fig:sensitivity-plots}, we observe poor performance for the majority of $\omega_1$ values when $\omega_0$ is ``too small''. However, for larger values of $\omega_0$, the performance becomes less sensitivity to perturbation in $\omega_1$. Multigrid solvers using inexact Braess-Sarazin relaxation have many more parameters than Vanka, which makes the sensitivity analysis slightly more complicated. For simplicity, we fix $\tau, {\alpha_i}$, and ${\beta_i}$ to be the same as in cycle BS$(2,2,1)$ in~\cref{table:LFA-mixed-two-grid1}, while varying $\omega_0$ and $\omega_1$. As seen in the right plot of~\cref{fig:sensitivity-plots}, there is a large region where we see near-optimal BS$(2,2,1)$ convergence without having to re-optimize the $\alpha_i$ and $\beta_i$. This suggests that the performance of BS$(2,2,1)$ is not too sensitive to its parameters.  We note that, for both relaxation schemes, there is a clear advantage to an informed choice of the relaxation parameters over naive choices, such as $\omega_0 = \omega_1 = 1$, but the sensitivity to their choices is not so severe as to be a practical drawback to using these approaches.

\subsection{Multilevel results}\label{sec:multilevel-results}

We now consider five-level algorithms, where the coarsest grid is a $4 \times 4$ mesh. The test problem is the same as in the previous subsection, with a uniform grid on the unit square and $h=\frac{1}{64}$.
We use the same parameters from~\cref{table:LFA-mixed-two-grid1} and compare the multilevel convergence results against both the LFA predicted two-grid convergence factors and the two-level performance reported above. The multilevel algorithm is formed by recursively extending the $h$-multigrid hierarchy in operator $G_1$. For completeness, \cref{table:LFA-mixed-multigrid} measures both V- and W-cycle convergence factors for stationary multigrid iterations, as well as the number of iterations required to reduce the norm of the absolute residual based on $K_0$ to $10^{-10}$, using both stationary GMG iterations and multigrid preconditioned FGMRES\@.
The measured convergence factors for W-cycles presented in~\cref{table:LFA-mixed-multigrid} agree with the LFA predicted results above. The V-cycle convergence rates, on the other hand, tend to underperform for both periodic and Dirichlet boundary conditions. Even though the W-cycle preconditioned FGMRES solver tends to converge in the fewest number of iterations,
the difference in iteration counts is not large, suggesting that V-cycle preconditioned FGMRES may be the more computationally efficient approach~\cite{Trottenberg2001}, particularly in a parallel setting. We note that the large \highlight{white}{orange!100}{discrepancy} between the number of iterations required for convergence using stationary V-cycles and V-cycle preconditioned FGMRES suggests that there are only a few modes that are not well-captured by the GMG operators.
\begin{table}[!ht]
  \caption{Predicted two-grid LFA convergence factors and measured convergence estimates for the multilevel method applied to the unit-square test problem.  A $P$ superscript denotes use of periodic boundary conditions, while no superscript denotes Dirichlet boundary conditions.  The last four columns list iteration counts to achieve a relative residual reduction tolerance using either stationary GMG iterations or GMG as a preconditioner for FGMRES.}
\centering
\begin{tabular}{llll c|c cc|cccc}
\toprule
&& & & & & & & \multicolumn{2}{c}{GMG} & \multicolumn{2}{c}{w/ FGMRES}\\
&$\nu_1$ & $\nu_2$ & $\gamma$& Cycle & $\hat{\rho}$ & $\rho^P$   & $\rho$ & $m^P$ & $m$ & $m^P$ & $m$   \\
\midrule
\multirow{4}{5mm}{BS} &\multirow{2}{*}{$1$} & \multirow{2}{*}{$0$} & \multirow{2}{*}{$1$} &V & \multirow{2}{*}{.18} &.17&.20&13&14&12&12\\
                      &  &&                         &W &                      &.17&.17&13&13&12&12\\  \cmidrule{2-12}
                      &\multirow{2}{*}{$1$} & \multirow{2}{*}{$1$} & \multirow{2}{*}{$1$} &V & \multirow{2}{*}{.09} &.25&.55&13&\highlight{white}{orange!100}{33}&9& \highlight{white}{orange!100}{11}\\
                      &  &&                         &W &                      &.08&.08&9&10&8&9\\
\midrule
\multirow{4}{5mm}{V}  &\multirow{2}{*}{$1$} & \multirow{2}{*}{$0$} & \multirow{2}{*}{$2$} &V & \multirow{2}{*}{.21} &.47&.20&21&14&12&11\\
                      &  &&                         &W &                      &.19&.20&14&14&11&11\\ \cmidrule{2-12}
                      &\multirow{2}{*}{$1$} & \multirow{2}{*}{$1$} & \multirow{2}{*}{$2$} &V & \multirow{2}{*}{.12} &.18&.18&12&12&10&10\\
                      &  &&                         &W &                      &.10&.09&10&10&9 &9 \\
\bottomrule
\end{tabular}\label{table:LFA-mixed-multigrid}
\end{table}

We note that these results use damping parameters obtained from the two-level analysis. As the two-level analysis assumes exact coarse-grid solves and W-cycles are expected to provide more accurate coarse solution than V-cycles, the W-cycle case is somewhat closer to the assumptions underlying the two-level analysis.  From this viewpoint, the better performance of the W-cycles is not surprising.  Given the observed success using V-cycles as preconditioners for FGMRES, we have not undertaken a multilevel LFA to see if there might possibly be better parameter choices for the V-cycle.

\subsection{Stokes backward-facing step problem}

In this section, we consider the performance of the preconditioned FGMRES algorithms as we increase the fine-scale problem size and number of levels in the multigrid hierarchy.  For this test, we consider a backward-facing step domain with the solution strategies discussed in~\cref{sec:multilevel-results} and verify that the resulting iteration counts are relatively stable with respect to mesh refinement. For all tests, the Stokes problem is defined on a backward-facing step domain with a parabolic inflow and natural outflow boundary conditions. A  Poiseuille  flow  profile  is imposed  on  the  inflow  boundary ($x=0; 0 \leq y \leq 2$). A homogeneous Neumann boundary condition is imposed on the outflow boundary ($x=2; 0 \leq y \leq 2$), thereby fixing the mean outflow pressure to be zero.
No-slip (zero velocity) boundary conditions are imposed on all other boundary faces.
More details on this type of problem can be found, for example, in~Elman et al.~\cite{elman2014finite}.

\Cref{table:refinement-study} summarizes the results, presenting iteration counts for multigrid preconditioned FGMRES\@. Using both IBSR and Vanka as relaxation within the preconditioner, we perform V- and W-cycle convergence tests with respect to problem size. The relaxation parameters are picked based on the most successful multilevel convergence results in~\cref{sec:multilevel-results}. The parameter values for these two multigrid hierarchies can be found in~\cref{table:LFA-mixed-two-grid1,table:LFA-mixed-two-grid2} under the corresponding $(\nu_1,\nu_2, \gamma)$ values:  $(1,1,1)$ for Braess-Sarazin and $(1,1,2)$ for Vanka relaxation.
For both IBSR and Vanka relaxation, the V-cycle results exhibit modest iteration growth as the problem size is increased. The convergence results for the W-cycle tests show constant iteration counts.
\begin{table}[!ht]
 \caption{FGMRES iterations to convergence for backward-facing step problem.
   Braess-Sarazin relaxation parameters are based on the results in~\cref{table:LFA-mixed-two-grid1}.
 Vanka relaxation parameters are based on the results in~\cref{table:LFA-mixed-multigrid}.}
\begin{center}
\begin{tabular}{r c c c c c}
\toprule
& &
\multicolumn{4}{c}{Method, $(\nu_1,\nu_2,\gamma)$}\\
& &
\multicolumn{2}{c}{BS, $(1,1,1)$} &
\multicolumn{2}{c}{V, $(1,1,2)$} \\
\cmidrule(lr){3-4} \cmidrule(lr){5-6}
DoFs         &GMG levels &V-cycle  &W-cycle  &V-cycle  &W-cycle \\
\midrule
\num{515   } &2          &10       &10       &9        &10      \\
\num{1891  } &3          &11       &10       &10       &10      \\
\num{7235  } &4          &12       &10       &11       &10      \\
\num{28291 } &5          &13       &10       &11       &10      \\
\num{111875} &6          &14       &10       &12       &10      \\
\num{444931} &7          &16       &10       &14       &10      \\
\bottomrule
\end{tabular}\label{table:refinement-study}
\end{center}
\end{table}

\begin{figure}[!ht]
  \centering
  \includegraphics{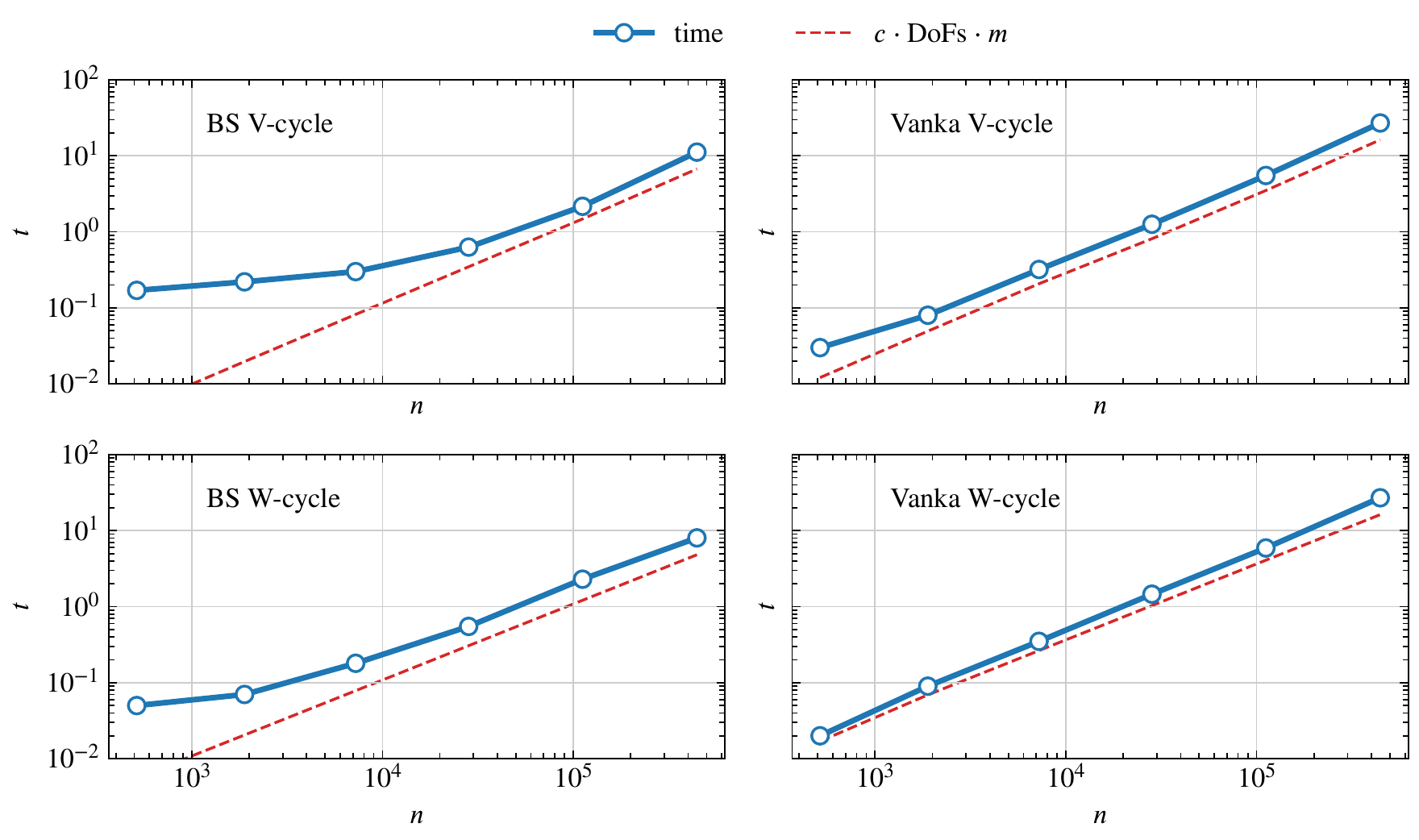}
  \caption{Timing of Braess-Sarazin and Vanka relaxation-based V- and W-cycles.  The dashed line represents the scaling with problem size in the form of $c \cdot \text{DoFs} \cdot m$, for some constant $c$ and where $m$ is the number of FGMRES iterations.}\label{fig:BVscaling}
\end{figure}

In addition to the number of FGMRES iterations in~\Cref{table:refinement-study}, we present the solve phase timings for Vanka and Braess-Sarazin based solvers in~\Cref{fig:BVscaling}. Both solvers's timings scale proportionally with the number of DoFs in the system and the total number of FGMRES iterations. While a single W-cycle is more computationally expensive per iteration than a V-cycle, the degradation in performance when using V-cycles makes W-cycles the faster option, due to the constant number of iterations to convergence. This is true for both Braess-Sarazin and Vanka based multigrid cycles.   While we see generally faster times for the Braess-Sarazin cycles here, we refrain from drawing conclusions about the relative performance between the two relaxation schemes, because timings are highly dependent on the specifics of their implementation and optimization~\cite{farrell2019pcpatch}.

\section{Conclusions and Future Work}\label{sec:conclusion}

In this paper, we demonstrate that the low-order, \isoqtwoqone{}, finite-element discretization can be used to construct effective preconditioners for the higher-order \qtwoqone{} finite-element discretization of the Stokes equations. To achieve effective performance, we use LFA in combination with robust optimization algorithms to compute relaxation parameters that optimize the resulting two-grid convergence factors for various two-grid algorithms based on the \isoqtwoqone{} hierarchy of grids.
The measured two-grid convergence factors for both periodic and Dirichlet boundary conditions generally agree quite closely with the LFA-predicted two-grid convergence factors. For multilevel convergence, we observe close agreement between the measured W-cycle convergence factors and the LFA predictions. The V-cycle multilevel convergence factors, however, can deviate significantly more from the LFA predictions than the W-cycle.  Both V- and W-cycles, however, lead to effective preconditioners for FGMRES for some parameter choices.  When used on a more challenging backward-facing step problem, W-cycle preconditioned FGMRES leads to no growth in iteration counts as the mesh is refined, while modest growth in iterations is seen with V-cycle preconditioners. Hence for the large problem sizes, W-cycle preconditioned FGMRES solvers are faster.

An immediate next step in this research is to see whether these results can be leveraged in the development of monolithic AMG algorithms for the \qtwoqone{} discretization of the Stokes equations.  Despite a long history of research effort, there has been little success in developing true algebraic MG approaches for saddle-point systems such as these.  Preliminary numerical results, that we intend to report in a future manuscript, suggest that it is easier to develop AMG algorithms for the \isoqtwoqone{} discretization as we consider here.  If successful, further work is possible for higher-order discretizations, such as the Scott-Vogelius discretization of the (Navier-) Stokes equations~\cite{LRScott_MVogelius_1985a, PEFarrell_etal_2021a}, or for coupled systems of equations of saddle-point type, such as viscoresistive magnetohydrodynamics~\cite{JAdler_etal_2020a}.

The source code used to collect data for this paper is publicly available at \url{https://github.com/lexeyV/Stokes_isoQ2Q1}. It is implemented in Python 3 under the 3-Clause BSD License. Version 1.0 is used in this paper, under commit \texttt{6ba1fa393dcdca113a011453630695afa49c9dfd}.

\section*{Acknowledgments}

The work of SPM was partially supported by an NSERC Discovery Grant. RT was supported by the U.S.~Department of Energy, Office of Science, Office of Advanced Scientific Computing Research, Applied Mathematics program.  Sandia National Laboratories is a multimission laboratory managed and operated by National Technology and Engineering Solutions of Sandia, LLC., a wholly owned subsidiary of Honeywell International, Inc., for the U.S. Department of Energy's National Nuclear Security Administration under grant DE-NA-0003525. This paper describes objective technical results and analysis. Any subjective views or opinions that might be expressed in the paper do not necessarily represent the views of the U.S. Department of Energy or the United States Government.

\bibliography{paper_isoq2q1_strip.bib}

\end{document}